# Some remarks on Rogers-Szegö polynomials and Losanitsch's triangle

Johann Cigler

Fakultät für Mathematik

Universität Wien

johann.cigler@univie.ac.at

**Abstract**

In this expository paper we collect some simple facts about analogues of Pascal's triangle where the entries count subsets of the integers with an even or odd sum and show that they are related to $q$ – Newton and Rogers-Szegö polynomials. In particular we consider an interesting triangle due to Losanitsch from this point of view. We also sketch some extensions of these results to sets whose sums have fixed residues modulo a prime $p$.

## 1. Introduction

This expository paper has been inspired by a Pascal-like triangle which has been obtained by the chemist S. Losanitsch (Sima Lozanic) [6]. Some of its properties can be found in [1] and [7] A034851. We show that the elements $L(n,k)$ of Losanitsch's triangle have several interesting interpretations:

a) For $k \equiv \pm 1 \bmod 4$ $L(n,k)$ is the number of $k$ – subsets of $\{1,2,\cdots,n\}$ with an even sum and for $k \not\equiv \pm 1 \bmod 4$ the number of $k$ – subsets of $\{1,2,\cdots,n\}$ with an odd sum.

b) Let $W_{n,k}$ be the set of all $n$ – tuples $w = (\varepsilon_1, \varepsilon_2, \cdots, \varepsilon_n)$ with $\varepsilon_i \in \{0,1\}$ and $\sum_{i=1}^{n} \varepsilon_i = k$. Then $L(n,k)$ is the number of those $n$ – tuples which have an even number of inversions.

c) For $w = (\varepsilon_1, \varepsilon_2, \cdots, \varepsilon_n) \in W_{n,k}$ let $\overline{w} = (\varepsilon_n, \cdots, \varepsilon_1)$ be the reversal of $w$. Then $L(n,k)$ also counts the set of all equivalence classes $\{w, \overline{w}\}$.

d) Let $\begin{bmatrix} n \\ k \end{bmatrix}$ be a $q$ – binomial coefficient. Then

$$\begin{bmatrix} n \\ k \end{bmatrix} \equiv \left( L(n,k) + \left( \binom{n}{k} - L(n,k) \right) q \right) \bmod (q^2 - 1).$$

For helpful remarks I want to thank J. Schoissengeier and V. Losert.

**1.1.**
As is well known the entries $\binom{n}{k}$ of Pascal's triangle count the subsets of $\{1,2,\cdots,n\}$ with $k$ elements which will be called $k$ – sets for short.
A $k$ – set $S$ will be called *even*, if the sum $\sigma(S)$ of its elements is even and *odd* if this sum is odd. By convention the empty set is even.

Let $E_{n,k}$ be the set of all even $k$ – subsets of $\{1,2,\cdots,n\}$ and $e(n,k) = |E_{n,k}|$ the number of its elements and let $O_{n,k}$ be the set of all odd $k$ – subsets of $\{1,2,\cdots,n\}$ and $o(n,k) = |O_{n,k}|$ the number of its elements.

For example $e(5,3) = 6$ because $E_{5,3} = \{\{1,2,3\},\{1,2,5\},\{1,3,4\},\{1,4,5\},\{2,3,5\},\{3,4,5\}\}$.

Let us note the trivial fact
$$e(n,k) + o(n,k) = \binom{n}{k}. \tag{1.1}$$

**Lemma 1.1**

$$e(n,k) = e(n-2,k) + \binom{n-2}{k-1} + o(n-2,k-2),$$
$$o(n,k) = o(n-2,k) + \binom{n-2}{k-1} + e(n-2,k-2). \tag{1.2}$$

This is a refinement of the formula $\binom{n}{k} = \binom{n-2}{k} + 2\binom{n-2}{k-1} + \binom{n-2}{k-2}$.

**Proof**

Let $S$ be a $k$ – subset of $\{1,2,\cdots,n\}$. To prove (1.2) consider 3 possibilities:

a) $S \subseteq \{1,2,\cdots,n-2\}$,

b) $S$ contains precisely one of the numbers $n-1$ and $n$. The remaining $(k-1)$ – set can be an arbitrary subset of $\{1,\cdots,n-2\}$ because one of these numbers is even and the other one is odd. There are $\binom{n-2}{k-1}$ such subsets.

c) $S$ contains both $n-1$ and $n$. Since $n+(n-1)$ is odd the remaining $(k-2)$ – subset must have the opposite parity of the given $k$ – subset.

**1.2**
To each $k$ – set $S$ we associate its indicator function $c_S$ defined by $c_S(i) = 1$ if $i \in S$ and $c_S(i) = 0$ else, or equivalently the $n$ – tuple $w(S) = (\varepsilon_1, \varepsilon_2, \cdots, \varepsilon_n)$ with $\varepsilon_i = c_S(i)$. . Let $W_{n,k}$ be the set of all $n$ – tuples $w = (\varepsilon_1, \varepsilon_2, \cdots, \varepsilon_n)$ with $\varepsilon_i \in \{0,1\}$ and $\sum_{i=1}^{n} \varepsilon_i = k$. For each $w$ let $\text{inv}(w) = |\{1 \le i < j \le n : \varepsilon_i = 1 > \varepsilon_j = 0\}|$ be the number of inversions of $w$. We call an $n$ – tuple even if $\text{inv}(w)$ is even and odd else.

For example there are 4 even and 2 odd quadruples $(\varepsilon_1, \varepsilon_2, \varepsilon_3, \varepsilon_4)$ with $\varepsilon_1 + \varepsilon_2 + \varepsilon_3 + \varepsilon_4 = 2$ because $\operatorname{inv}((0,0,1,1)) = 0$, $\operatorname{inv}((0,1,0,1)) = 1$, $\operatorname{inv}((0,1,1,0)) = 2$, $\operatorname{inv}((1,0,0,1)) = 2$, $\operatorname{inv}((1,0,1,0)) = 3$, $\operatorname{inv}((1,1,0,0)) = 4$.

Instead of $W_{n,k}$ we could also consider the set of all products $w = C_1 C_2 \cdots C_n$ of $k$ elements $A$ and $n-k$ elements $B$ which satisfy $BA = qAB$. An inversion is then a pair $(C_i, C_j) = (B, A)$ with $i < j$ which gives $C_1 C_2 \cdots C_n = q^{\operatorname{inv}(w)} A^k B^{n-k}$.

Let us associate to each $S = \{i_1, i_2, \cdots, i_k\}$ with $1 \leq i_1 < i_2 < \cdots < i_k \leq n$ the set
$S^* = \{n+1-i_k, n+1-i_{k-1}, \cdots, n+1-i_1\}$.

**Lemma 1.2**

For each $k-$set $S \subseteq \{1, 2, \cdots, n\}$ we have $\sigma(S) = \binom{k+1}{2} + \operatorname{inv}(w(S^*))$.

**Proof**

This holds for $k = 1$ and all $n$ because $S = \{i\}$ implies $S^* = \{n+1-i\}$ and $\sigma(S) = i$ and $\operatorname{inv}(w(S^*)) = i - 1$. Thus $i = \binom{1+1}{2} + i - 1$.

Let us assume that it is true for $k-1$ and all $n$ and consider a $k-$set $S$.

It is trivially true for $n = k$. Let us suppose that it is true for $n-1$.

Let $S_0$ be the restriction of $S$ to $\{1, 2, \cdots, n-1\}$. If $S_0^* = (\varepsilon_1, \cdots, \varepsilon_{n-1})$ then $S^* = \{c_S(n), \varepsilon_1, \cdots, \varepsilon_{n-1}\}$.

If $n \notin S$ then $\sigma(S) = \sigma(S_0) = \binom{k+1}{2} + \operatorname{inv}(w(S_0^*)) = \binom{k+1}{2} + \operatorname{inv}(w(S^*))$.

If $n \in S$ then
$\sigma(S) = \sigma(S_0) + n = \binom{k}{2} + \operatorname{inv}(w(S_0^*)) + n = \binom{k+1}{2} + n - k + \operatorname{inv}(w(S_0^*)) = \binom{k+1}{2} + \operatorname{inv}(w(S^*))$
because there are $n-k$ elements $0$ in $w(S_0^*)$.

We will now consider matrices whose columns are either $e_k$ or $o_k$, where $e_k$ is the column with entries $e(n,k)$ and $o_k$ the column with entries $o(n,k)$ for $n \in \mathbb{N}$.

## 2. Sets with an even sum

### 2.1.

Let us first consider the matrix $(e(n,k)) = (e_0, e_1, e_2, e_3, \cdots)$ whose entries are the number of even sets (cf. OEIS [7], A282011). The first terms are

$$\begin{pmatrix} 1 & 0 & 0 & 0 & 0 & 0 & 0 \\ 1 & 0 & 0 & 0 & 0 & 0 & 0 \\ 1 & 1 & 0 & 0 & 0 & 0 & 0 \\ 1 & 1 & 1 & 1 & 0 & 0 & 0 \\ 1 & 2 & 2 & 2 & 1 & 0 & 0 \\ 1 & 2 & 4 & 6 & 3 & 0 & 0 \\ 1 & 3 & 6 & 10 & 9 & 3 & 0 \end{pmatrix}$$

Closely related is the matrix $(o(n,k)) = (o_0, o_1, o_2, o_3, \cdots)$ whose entries are the number of odd sets (cf. OEIS [7], A 159916).

$$\begin{pmatrix} 0 & 0 & 0 & 0 & 0 & 0 & 0 \\ 0 & 1 & 0 & 0 & 0 & 0 & 0 \\ 0 & 1 & 1 & 0 & 0 & 0 & 0 \\ 0 & 2 & 2 & 0 & 0 & 0 & 0 \\ 0 & 2 & 4 & 2 & 0 & 0 & 0 \\ 0 & 3 & 6 & 4 & 2 & 1 & 0 \\ 0 & 3 & 9 & 10 & 6 & 3 & 1 \end{pmatrix}$$

Since in both triangles the columns $c_k$ and $c_{k+2}$ have the same parity, by Lemma 1.1 the entries $a(n,k)$ of these matrices satisfy

$$a(n,k) = a(n-2, k) + \binom{n-2}{k-1} + \binom{n-2}{k-2} - a(n-2, k-2) \text{ and thus}$$

$$a(n,k) = a(n-2, k) + \binom{n-1}{k-1} - a(n-2, k-2). \tag{2.1}$$

Therefore the polynomials $a_n(x) = \sum_{k=0}^{n} a(n,k) x^k$ satisfy the recursion

$$a_n(x) = (1 - x^2) a_{n-2}(x) + x(1+x)^{n-1}. \tag{2.2}$$

By applying this to $n$ and $n-1$ we get the homogeneous recursion

$$a_n(x) = (1+x) a_{n-1}(x) + (1 - x^2) a_{n-2}(x) - (1+x)(1 - x^2) a_{n-3}(x). \tag{2.3}$$

Since $z^3 - (1+x) z^2 - (1 - x^2) z + (1+x)(1 - x^2) = (z - (1+x))\left(z - \sqrt{1 - x^2}\right)\left(z + \sqrt{1 - x^2}\right)$

we see that $a_n(x) = c_0 (1+x)^n + c_1 \left(\sqrt{1 - x^2}\right)^n + c_2 \left(-\sqrt{1 - x^2}\right)^n$ for some coefficients $c_i$.

Let now $e_n(x) = \sum_{k=0}^{n} e(n,k) x^k$ and $o_n(x) = \sum_{k=0}^{n} o(n,k) x^k$.

From the initial values $e_0(x) = 1, e_1(x) = 1, e_2(x) = 1+x$ we compute

$$c_0 = \frac{1}{2}, c_1 = \frac{1-x+\sqrt{1-x^2}}{4\sqrt{1-x^2}}, c_2 = \frac{-1+x+\sqrt{1-x^2}}{4\sqrt{1-x^2}}.$$

This gives

$$e_{2n}(x) = \frac{(1+x)^{2n} + (1-x^2)^n}{2},$$
$$e_{2n+1}(x) = \frac{(1+x)^{2n+1} + (1-x)(1-x^2)^n}{2} \qquad (2.4)$$

and analogous formulae for $o_n(x)$.

From (2.4) we see that $e_n(x)$ is palindromic, i.e. $x^n e_n\left(\frac{1}{x}\right) = e_n(x)$ if $n \equiv 0, 3 \bmod 4$.

The generating functions are

$$\sum_{n\geq 0} e_n(x) z^n = \frac{1 - xz - (1-x^2)z^2}{(1-(1+x)z)(1-(1-x^2)z^2)} = \frac{1}{2}\left(\frac{1}{1-(1+x)z} + \frac{1+(1-x)z}{1-(1-x^2)z^2}\right) \qquad (2.5)$$

and

$$\sum_{n\geq 0} o_n(x) z^n = \frac{xz}{(1-(1+x)z)(1-(1-x^2)z^2)} = \frac{1}{2}\left(\frac{1}{1-(1+x)z} - \frac{1+(1-x)z}{1-(1-x^2)z^2}\right) \qquad (2.6)$$

From (2.5) we conclude that the generating function of $e_n^*(x) = \sum_{k=0}^{n} e(n, n-k) x^k = x^n e_n\left(\frac{1}{x}\right)$ is

$$\sum_{n\geq 0} e_n^*(x) z^n = \frac{1 - z - (x^2-1)z^2}{(1-(1+x)z)(1+(1-x^2)z^2)} = \frac{1}{2}\left(\frac{1}{1-(1+x)z} + \frac{1+(x-1)z}{1-(x^2-1)z^2}\right). \qquad (2.7)$$

**2.2**
It is now easy to derive some explicit formulae. From (2.4) we see that

$$e_n(x) = (1+x)^{\lfloor \frac{n}{2} \rfloor} \frac{(1+x)^{\lfloor \frac{n+1}{2} \rfloor} + (1-x)^{\lfloor \frac{n+1}{2} \rfloor}}{2} = (1+x)^{\lfloor \frac{n}{2} \rfloor} \sum_j \binom{\lfloor \frac{n+1}{2} \rfloor}{2j} x^{2j}. \qquad (2.8)$$

This implies

$$e(n,k) = \sum_j \binom{\left\lfloor \frac{n+1}{2} \right\rfloor}{2j} \binom{\left\lfloor \frac{n}{2} \right\rfloor}{k-2j}. \tag{2.9}$$

As special case we get the well-known formula

$$e(2n,n) = \sum_k \binom{n}{2k}^2. \tag{2.10}$$

If we write (2.5) in the form

$$\sum_{n \geq 0} \sum_k e(n,k) x^k z^n = \frac{1}{2} \left( \frac{1}{1-z} \frac{1}{1 - \frac{xz}{1-z}} + \frac{1}{1-z^2} \frac{1+z-xz}{1 + \frac{x^2 z^2}{1-z^2}} \right)$$

we see by comparing coefficients of $x$ that

$$\sum_{n \geq 0} e(n, 2k) z^n = \frac{z^{2k}}{2} \left( \frac{1}{(1-z)^{2k+1}} + \frac{(-1)^k (1+z)}{\left(1-z^2\right)^{k+1}} \right) \tag{2.11}$$

$$\sum_{n \geq 0} e(n, 2k+1) z^n = \frac{z^{2k+1}}{2} \left( \frac{1}{(1-z)^{2k+2}} + \frac{(-1)^{k+1}}{\left(1-z^2\right)^{k+1}} \right) \tag{2.12}$$

**Remark**
It is perhaps interesting that these generating functions can also be written in the following way:

$$\begin{aligned}
\sum_{n \geq 0} e(n, 4k) x^n &= \frac{x^{4k} e_{4k}(x)}{(1+x)^{4k}(1-x)^{4k+1}}, \\
\sum_{n \geq 0} e(n, 4k+1) x^n &= \frac{x^{4k+1} o_{4k+1}(x)}{(1+x)^{4k+1}(1-x)^{4k+2}}, \\
\sum_{n \geq 0} e(n, 4k+2) x^n &= \frac{x^{4k+2} o_{4k+2}(x)}{(1+x)^{4k+2}(1-x)^{4k+3}}, \\
\sum_{n \geq 0} e(n, 4k+3) x^n &= \frac{x^{4k+3} e_{4k+3}(x)}{(1+x)^{4k+3}(1-x)^{4k+4}}.
\end{aligned} \tag{2.13}$$

## 3. Losanitsch's triangle

### 3.1

Let us now consider another class of triangles where the columns $c_k$ and $c_{k+2}$ have opposite parity.

By Lemma 1.1 the entries of these matrices satisfy

$$b(n,k) = b(n-2,k) + \binom{n-2}{k-1} + b(n-2, k-2). \tag{3.1}$$

Therefore the polynomials $b_n(x) = \sum_{k=0}^{n} b(n,k) x^k$ satisfy the recursion

$$b_n(x) = (1+x^2) b_{n-2}(x) + x(1+x)^{n-2}. \tag{3.2}$$

By applying this to $n$ and $n-1$ we get

$$b_n(x) = (1+x) b_{n-1}(x) + (1+x^2) b_{n-2}(x) - (1+x)(1+x^2) b_{n-3}(x). \tag{3.3}$$

The best known special case of a matrix satisfying (3.1) is **Losanitsch's triangle** $(L(n,k))$ (cf. [1], [7] A034851) which is defined by the recursion

$$L(n,k) = L(n-2,k) + \binom{n-2}{k-1} + L(n-2, k-2) \tag{3.4}$$

with initial values $L(0,k) = [k=0]$ and $L(1,k) = [k \leq 1]$ and $L(n,k) = 0$ for $k < 0$.

The first terms are

$$\begin{pmatrix} 1 & 0 & 0 & 0 & 0 & 0 & 0 \\ 1 & 1 & 0 & 0 & 0 & 0 & 0 \\ 1 & 1 & 1 & 0 & 0 & 0 & 0 \\ 1 & 2 & 2 & 1 & 0 & 0 & 0 \\ 1 & 2 & 4 & 2 & 1 & 0 & 0 \\ 1 & 3 & 6 & 6 & 3 & 1 & 0 \\ 1 & 3 & 9 & 10 & 9 & 3 & 1 \end{pmatrix}$$

This matrix has been obtained by the chemist S.M. Losanitsch (or Lozanic) [6] in his investigation of paraffin. Therefore we call the numbers $L(n,k)$ Losanitsch numbers. The same triangle has also been considered in [1] in the study of some sort of necklaces where these numbers have been called necklace numbers. Further information can be found in OEIS [7], A034851.

Let us first observe that

$$(L(n,k)) = (e_0, o_1, o_2, e_3, e_4, o_5, \cdots). \tag{3.5}$$

For the right-hand side satisfies (3.4) with the same initial and boundary values as Losanitsch's triangle.

By (2.4) we see that $e(n,n)=1$ if $n \equiv 0,3 \mod 4$ and $o(n,n)=1$ if $n \equiv 1,2 \mod 4$. Therefore we get

**Proposition 3.1**
*Losanitsch's triangle is the uniquely determined matrix where for each $k$ the $k-$th column is either $e_k$ or $o_k$ and all elements of the main diagonal are 1.*

Let us give another characterization of Losanitsch's triangle which appears as special case in [5].

**Theorem 3.2** (Stephen G. Hartke and A.J. Radcliffe)

*The Losanitsch number $L(n,k)$ is the number of all $n-$tuples $L_{n,k} \subseteq W_{n,k}$ with an even number of inversions.*

**Proof**

By Lemma 1.2 we see that the map $S \to w(S^*)$ from all $k-$subsets of $\{1,\cdots,n\}$ to $W_{n,k}$ is a bijection between the $k-$sets for which $\sigma(S) \equiv \binom{k+1}{2} \mod 2$ and $L_{n,k}$.

For $k \equiv \pm 1 \mod 4$ we have $\binom{k+1}{2} \equiv 0 \mod 2$ and therefore $|L_{n,k}| = e(n,k)$ and for $k \equiv 1,2 \mod 4$ we have $\binom{k+1}{2} \equiv 1 \mod 2$ and therefore $|L_{n,k}| = o(n,k)$.

Thus $|L_{n,k}| = L(n,k)$.

**3.2**
The opposite matrix $(\bar{L}(n,k)) = (o_0, e_1, e_2, o_3, o_4, e_5, \cdots)$ is OEIS [7], A034852 and essentially also A034877. The first terms are

$$\begin{pmatrix} 0 & 0 & 0 & 0 & 0 & 0 & 0 \\ 0 & 0 & 0 & 0 & 0 & 0 & 0 \\ 0 & 1 & 0 & 0 & 0 & 0 & 0 \\ 0 & 1 & 1 & 0 & 0 & 0 & 0 \\ 0 & 2 & 2 & 2 & 0 & 0 & 0 \\ 0 & 2 & 4 & 4 & 2 & 0 & 0 \\ 0 & 3 & 6 & 10 & 6 & 3 & 0 \end{pmatrix}$$

The Losanitsch polynomials $L_n(x) = \sum_{k=0}^{n} L(n,k)x^k$ and $\bar{L}_n(x) = \sum_{k=0}^{n} \bar{L}(n,k)x^k$ satisfy the recursion (3.3).

In the same way as above we see that

$$L_n(x) = \frac{(1+x)^n}{2} + \frac{1+x+\sqrt{1+x^2}}{4\sqrt{1+x^2}}\left(\sqrt{1+x^2}\right)^n + \frac{-1-x+\sqrt{1+x^2}}{4\sqrt{1+x^2}}\left(-\sqrt{1+x^2}\right)^n.$$

This can be simplified to

$$L_{2n}(x) = \frac{(1+x)^{2n} + \left(1+x^2\right)^n}{2},$$

$$L_{2n+1}(x) = \frac{(1+x)^{2n+1} + (1+x)\left(1+x^2\right)^n}{2}.$$

(3.6)

Note that

$$L_{2n+1}(x) = (1+x)L_{2n}(x) \tag{3.7}$$

and

$$x^n L_n\left(\frac{1}{x}\right) = L_n(x) \tag{3.8}$$

or

$$L(n, n-k) = L(n, k). \tag{3.9}$$

Thus we get (cf. [1], Theorem 2.8) by observing (1.1)

$$L(2n, 2k+1) = \bar{L}(2n, 2k+1) = \frac{1}{2}\binom{2n}{2k+1},$$

$$L(n,k) = \frac{1}{2}\left(\binom{n}{k} + \binom{\lfloor\frac{n}{2}\rfloor}{\lfloor\frac{k}{2}\rfloor}\right) \text{ and } \bar{L}(n,k) = \frac{1}{2}\left(\binom{n}{k} - \binom{\lfloor\frac{n}{2}\rfloor}{\lfloor\frac{k}{2}\rfloor}\right) \text{ else.}$$

(3.10)

In the same way as above we get the generating function

$$\sum_{n\geq 0} L_n(x)z^n = \frac{1-\left(1+x+x^2\right)z^2}{\left(1-(1+x)z\right)\left(1-\left(1+x^2\right)z^2\right)} = \frac{1}{2}\left(\frac{1}{1-(1+x)z} + \frac{1+(1+x)z}{1-\left(1+x^2\right)z^2}\right). \tag{3.11}$$

Some other properties of the Losanitsch polynomials can be found in [1].

Comparing with (2.13) we get

**Proposition 3.3**

$$\sum_n L(n,k)x^n = \sum_n L(n,n-k)x^n = \frac{x^k e_k(x)}{(1-x)^{k+1}(1+x)^k}. \tag{3.12}$$

There exists also another interesting relation between the numbers $e(n,k)$ and $L(n,k)$.

**Proposition 3.4**

$$\sum_{n} e(n, n-k) z^n = \frac{z^k}{(1-z)^{2\left\lfloor\frac{k+1}{2}\right\rfloor+1} \left(1+z^2\right)^{\left\lfloor\frac{k}{2}\right\rfloor+1}} L_{k+2}(-z). \tag{3.13}$$

**Proof**

It suffices to compute the coefficient of $x^k$ in (2.7). This gives

$$\sum_{n \geq 0} e(n, n-2k) z^n = \frac{1}{2}\left( \frac{z^{2k}}{(1-z)^{2k+1}} + \frac{z^{2k}(1-z)}{(1+z^2)^{k+1}} \right) = \frac{z^{2k}}{(1-z)^{2k+1}(1+z^2)^{k+1}} \left( \frac{(1+z^2)^{k+1} + (1-z)^{2k+2}}{2} \right)$$

$$= \frac{z^{2k}}{(1-z)^{2k+1}(1+z^2)^{k+1}} L_{2k+2}(-z)$$

and

$$\sum_{n \geq 0} e(n, n-2k-1) z^n = \frac{1}{2}\left( \frac{z^{2k+1}}{(1-z)^{2k+2}} + \frac{z^{2k+1}}{(1+z^2)^{k+1}} \right) = \frac{z^{2k+1}}{(1-z)^{2k+3}(1+z^2)^{k+1}} \left( \frac{(1-z)(1+z^2)^{k+1} + (1-z)^{2k+3}}{2} \right)$$

$$= \frac{z^{2k+1}}{(1-z)^{2k+3}(1+z^2)^{k+1}} L_{2k+3}(-z).$$

**3.3**
Let us now derive another interesting combinatorial interpretation of the Losanitsch numbers. For $w = (\varepsilon_1, \cdots, \varepsilon_n) \in W_{n,k}$ let $\overline{w} = (\varepsilon_n, \cdots, \varepsilon_1)$ be the reversal of $w$. Define an equivalence relation on $W_{n,k}$ by $w \sim \overline{w}$. If $w = \overline{w}$ we call $w$ palindromic.

**Theorem 3.5**

*Let $R_{n,k}$ be the set of all equivalence classes $\{w, \overline{w}\}$ and let $P_{n,k}$ be the subset of all palindromic classes. Let $r(n,k) = |R_{n,k}|$ and $p(n,k) = |P_{n,k}|$ denote the number of their elements. Then*

$$\begin{aligned} r(n,k) &= L(n,k), \\ p(n,k) &= L(n,k) - \overline{L}(n,k). \end{aligned} \tag{3.14}$$

**Proof**

By (3.10) we have $L(n, 2k+1) - \overline{L}(n, 2k+1) = 0$ and $L(n,k) - \overline{L}(n,k) = \left( \begin{array}{c} \lfloor \frac{n}{2} \rfloor \\ \lfloor \frac{k}{2} \rfloor \end{array} \right)$ else.

For $p(n,k)$ we get the same result: $p(2n, 2k) = \binom{n}{k}$ because each palindromic set is uniquely determined by its restriction to $\{1, 2, \cdots, n\}$ and each $k$-subset of $\{1, 2, \cdots, n\}$ has a unique extension to a palindromic set of $\{1, 2, \cdots, 2n\}$.

It is clear that $p(2n, 2k+1) = 0$.

$p(2n+1, 2k+1) = \binom{n}{k}$ because each $k$-subset of $\{1, 2, \cdots, n+1\}$ has a unique extension to a palindromic set of $\{1, 2, \cdots, 2n+1\}$.

Thus $p(n,k) = L(n,k) - \overline{L}(n,k)$.

Since there are $\frac{1}{2}\left( \binom{n}{k} - p(n,k) \right)$ non-palindromic classes in $W_{n,k}$ we get

$r(n,k) = p(n,k) + \frac{1}{2}\left( \binom{n}{k} - p(n,k) \right) = \frac{1}{2}\left( \binom{n}{k} + p(n,k) \right) = L(n,k)$.

Therefore $\overline{L}(n,k) = \frac{1}{2}\left( \binom{n}{k} - p(n,k) \right)$ is the number of non-palindromic equivalence classes.

**Example**

For example $W_{5,2}$ has 6 equivalence classes

$\{00011, 11000\}, \{00101, 10100\}, \{01001, 10010\}, \{01100, 00110\},$
$\{10001\}, \{01010\}$

Thus $L(5,2) = 6$ and $\overline{L}(5,2) = 4$ and $p(5,2) = L(5,2] - \overline{L}(5,2) = 2$.

**Remark**

Some special cases of Theorem 3.5 are mentioned in OEIS A002620, A005993, A005994 and A005995 in the following formulation:
Let $B(n,k)$ be the number of bracelets, i.e. necklaces allowing turning over, with $n$ white, 1 blue and $k$ red beads, then $B(n,k) = L(n+k, k)$.
To verify this it suffices to choose a representative of the bracelet whose first term is the blue bead.
If we replace the first $k-1$ letters 1 by 10 we see that Theorem 3.5 implies [1], Theorem 2.4.

## 4. Rogers-Szegö polynomials and Losanitsch's triangle

**4.1**

There are close connections with some $q-$binomial theorems. For the relevant facts about $q-$calculus I refer to [2]. Let $\begin{bmatrix} n \\ k \end{bmatrix} = \begin{bmatrix} n \\ k \end{bmatrix}_q$ be a $q-$binomial coefficient. It turns out that both $e(n,k)$ and $L(n,k)$ are residues modulo $q^2-1$ of $q-$binomial coefficients.

**Theorem 4.1**

$$e(n,k) + o(n,k)q \equiv q^{\binom{k+1}{2}} \begin{bmatrix} n \\ k \end{bmatrix} \mod(q^2-1),$$

$$L(n,k) + \bar{L}(n,k)q \equiv \begin{bmatrix} n \\ k \end{bmatrix} \mod(q^2-1).$$

(4.1)

**Proof**

Recall that

$$\prod_{j=1}^{n}(1+q^j x) = \sum_{k=0}^{n} q^{\binom{k+1}{2}} \begin{bmatrix} n \\ k \end{bmatrix} x^k.$$

(4.2)

Since $\prod_{j=1}^{n}(1+q^j x) = \sum_{k} \sum_{j_1 < \cdots < j_k} q^{j_1+j_2+\cdots+j_k} x^k \equiv \sum_{k}(e(n,k)+o(n,k)q) x^k \mod(q^2-1)$

we get the first line of (4.1).

Now observe that $\binom{4k+1}{2} \equiv \binom{4k+4}{2} \equiv 0 \mod 2$ and $\binom{4k+2}{2} \equiv \binom{4k+3}{2} \equiv 1 \mod 2$.

Therefore

$$\begin{bmatrix} n \\ k \end{bmatrix} \mod(q^2-1) \equiv q^{\binom{k+1}{2}} \begin{bmatrix} n \\ k \end{bmatrix} \mod(q^2-1) = e(n,k)+o(n,k)q \text{ for } k \equiv 0,3 \mod 4$$

and

$$\begin{bmatrix} n \\ k \end{bmatrix} \mod(q^2-1) \equiv qq^{\binom{k+1}{2}} \begin{bmatrix} n \\ k \end{bmatrix} \mod(q^2-1) = o(n,k)+qe(n,k) \text{ for } k \equiv 1,2 \mod 4.$$

Comparing with (3.5) we get the second line of (4.1).

**Corollary 4.2**

*The Losanitsch polynomials are related to the Rogers-Szegö polynomials* $r_n(x,q) = \sum_{k=0}^{n} \begin{bmatrix} n \\ k \end{bmatrix} x^k$

*by*

$$L_n(x) + q\overline{L}_n(x) \equiv r_n(x,q) \bmod (q^2 - 1) \quad (4.3)$$

*and the polynomials* $e_n(x)$ *are related to the* $q$ – *Newton polynomials*

$$p_n(x,q) = \prod_{j=1}^{n} (1 + q^j x) = \sum_{k=0}^{n} q^{\binom{k+1}{2}} \begin{bmatrix} n \\ k \end{bmatrix} x^k \text{ by}$$

$$e_n(x) + qo_n(x) \equiv p_n(x,q) \bmod (q^2 - 1). \quad (4.4)$$

From (4.3) it is obvious that $L_n(x)$ is palindromic.

We already know that $e_n(x)$ is palindromic if $n \equiv 0, 3 \bmod 4$.
This can also be deduced from (4.1).

For $n \equiv 3 \bmod 4$ this follows from $\binom{4n+4-k}{2} \equiv \binom{k+1}{2} \bmod 2$.

For $n \equiv 0 \bmod 4$ we have $\binom{4n+1-k}{2} \equiv k \bmod 2$. This implies that

$$q^{\binom{4n+1-k}{2}} \begin{bmatrix} 4n \\ 4n-k \end{bmatrix} \equiv q \begin{bmatrix} 4n \\ k \end{bmatrix} \equiv o(4n,k) + qe(4n,k) \bmod (q^2 - 1) \text{ if } k \equiv 1 \bmod 2. \text{ But in this case}$$

we have $e(4n,k) = o(4n,k)$ since by (2.4) $e_{2n}(x) - o_{2n}(x) = (1 - x^2)^n$.

**4.2**

Several theorems about $q$ – binomial coefficients give results for Losanitsch numbers or for the numbers of even or odd sets if we replace the polynomials $e_n(x)$ and $L_n(x)$ by the equivalence classes $\varepsilon_n(x) = e_n(x) + qo_n(x)$ and $\lambda_n(x) = L_n(x) + q\overline{L}_n(x)$ modulo $q^2 - 1$.

The first terms of the matrix $(e(n,k) + qo(n,k))$ are

$$\begin{pmatrix} 1 & 0 & 0 & 0 & 0 & 0 & 0 & 0 \\ 1 & q & 0 & 0 & 0 & 0 & 0 & 0 \\ 1 & 1+q & q & 0 & 0 & 0 & 0 & 0 \\ 1 & 1+2q & 1+2q & 1 & 0 & 0 & 0 & 0 \\ 1 & 2+2q & 2+4q & 2+2q & 1 & 0 & 0 & 0 \\ 1 & 2+3q & 4+6q & 6+4q & 3+2q & q & 0 & 0 \\ 1 & 3+3q & 6+9q & 10+10q & 9+6q & 3+3q & q & 0 \\ 1 & 3+4q & 9+12q & 19+16q & 19+16q & 9+12q & 3+4q & 1 \end{pmatrix}$$

Let us note that

$$(1-q)\varepsilon_n(x) \equiv (1-q)(1-x)^{\left\lfloor \frac{n+1}{2} \right\rfloor}(1+x)^{\left\lfloor \frac{n}{2} \right\rfloor} \mod(q^2-1). \tag{4.5}$$

This holds for $n=0$ and $n=1$. By induction we get

$$(1-q)\varepsilon_{2n}(x) \equiv (1-q)\prod_{j=1}^{2n}(1+q^j x) \equiv (1+q^{2n}x)(1-q)\varepsilon_{2n-1}(x)$$

$$\equiv (1+x)(1-q)(1-x)^n(1+x)^{n-1} \equiv (1-q)(1-x)^n(1+x)^n$$

and

$$(1-q)\varepsilon_{2n+1}(x) \equiv (1-q)\prod_{j=1}^{2n+1}(1+q^j x) \equiv (1+q^{2n+1}x)(1-q)\varepsilon_{2n}(x)$$

$$\equiv (1+qx)(1-q)(1-x)^n(1+x)^n \equiv (1-q)(1-x)(1-x)^n(1+x)^n \equiv (1-q)(1-x)^{n+1}(1+x)^n \mod(q^2-1).$$

Here we have used that $(1+qx)(1-q) \equiv 1+qx-q-x \equiv (1-x)(1-q) \mod(q^2-1)$.

The first terms of $\left(L(n,k)+q\overline{L}(n,k)\right)$ are

$$\begin{pmatrix}
1 & 0 & 0 & 0 & 0 & 0 & 0 & 0 \\
1 & 1 & 0 & 0 & 0 & 0 & 0 & 0 \\
1 & 1+q & 1 & 0 & 0 & 0 & 0 & 0 \\
1 & 2+q & 2+q & 1 & 0 & 0 & 0 & 0 \\
1 & 2+2q & 4+2q & 2+2q & 1 & 0 & 0 & 0 \\
1 & 3+2q & 6+4q & 6+4q & 3+2q & 1 & 0 & 0 \\
1 & 3+3q & 9+6q & 10+10q & 9+6q & 3+3q & 1 & 0 \\
1 & 4+3q & 12+9q & 19+16q & 19+16q & 12+9q & 4+3q & 1
\end{pmatrix}$$

In this case we get

$$(1-q)\lambda_{2n}(x) \equiv (1-q)(1+x^2)^n \mod(q^2-1),$$
$$\lambda_{2n+1}(x) \equiv (1+x)\lambda_{2n}(x). \tag{4.6}$$

The recurrence

$$r_n(x,q) = (x+1)r_{n-1}(x,q) + (q^{n-1}-1)xr_{n-2}(x,q) \tag{4.7}$$

for the Rogers-Szegö polynomials (cf. e.g. [2]) implies $\lambda_{2n+1}(x) \equiv (1+x)\lambda_{2n}(x)$ and
$\lambda_{2n}(x) \equiv (x+1)\lambda_{2n-1}(x) + (q^{2n-1}-1)x\lambda_{2n-2}(x)$. Since
$(1-q)(q^{2n-1}-1) \equiv -(1-q)^2 \equiv -2(1-q) \mod(q^2-1)$ we get
$(1-q)\lambda_{2n}(x) \equiv (x+1)(1-q)\lambda_{2n-1}(x) - 2(1-q)x\lambda_{2n-2}(x)$
$\equiv (x+1)^2(1-q)\lambda_{2n-2}(x) - 2(1-q)x\lambda_{2n-2}(x) \equiv (1-q)(x^2+1)\lambda_{2n-2}(x)$
which implies the first line of (4.6) by induction.

Now it is well known that

$$\begin{bmatrix} n \\ k \end{bmatrix}_q = \sum_{w \in W_{n,k}} q^{\mathrm{inv}(w)}. \qquad (4.8)$$

This follows for example from the $q-$binomial theorem $(A+B)^n = \sum_{k=0}^n \begin{bmatrix} n \\ k \end{bmatrix} A^k B^{n-k}$ for $q-$commuting variables $BA = qAB$. (cf. e.g. [2]). For if we write $(A+B)^n = \sum C_{i_1} \cdots C_{i_n}$ with $C_i \in \{A, B\}$ we get $(A+B)^n = \sum_k A^k B^{n-k} \sum_{w \in W_{n,k}} q^{\mathrm{inv}(w)}$.

If we reduce (4.8) modulo $q^2 - 1$ we get again Theorem 3.2.

**4.3**

Since each equivalence class modulo $q^2 - 1 = (q-1)(q+1)$ is also an equivalence class modulo $q-1$ and modulo $q+1$ we see that

$$L(n,k) + \overline{L}(n,k) = \begin{bmatrix} n \\ k \end{bmatrix}_{q=1} = \binom{n}{k} \text{ and } L(n,k) - \overline{L}(n,k) = \begin{bmatrix} n \\ k \end{bmatrix}_{q=-1}.$$

Thus

$$L(n,k) = \frac{1}{2}\left( \begin{bmatrix} n \\ k \end{bmatrix}_{q=1} + \begin{bmatrix} n \\ k \end{bmatrix}_{q=-1} \right). \qquad (4.9)$$

By (3.10) we get the well-known result

$$\begin{bmatrix} 2n \\ 2k+1 \end{bmatrix}_{q=-1} = 0,$$

$$\begin{bmatrix} n \\ k \end{bmatrix}_{q=-1} = \binom{\lfloor \frac{n}{2} \rfloor}{\lfloor \frac{k}{2} \rfloor} \text{ else.} \qquad (4.10)$$

The first terms of the table $\left( \begin{bmatrix} n \\ k \end{bmatrix}_{q=-1} \right)_{n,k \geq 0}$ are

$$\begin{pmatrix} 1 & 0 & 0 & 0 & 0 & 0 & 0 & 0 \\ 1 & 1 & 0 & 0 & 0 & 0 & 0 & 0 \\ 1 & 0 & 1 & 0 & 0 & 0 & 0 & 0 \\ 1 & 1 & 1 & 1 & 0 & 0 & 0 & 0 \\ 1 & 0 & 2 & 0 & 1 & 0 & 0 & 0 \\ 1 & 1 & 2 & 2 & 1 & 1 & 0 & 0 \\ 1 & 0 & 3 & 0 & 3 & 0 & 1 & 0 \\ 1 & 1 & 3 & 3 & 3 & 3 & 1 & 1 \end{pmatrix}$$

**Remark**

If we compare (4.10) with Theorem 3.5 we see that the number of palindromic $n-$ tuples in $W_{n,k}$ coincides with

$$p(n,k) = \begin{bmatrix} n \\ k \end{bmatrix}_{q=-1}. \tag{4.11}$$

(4.7) implies

$r_n(x,-1) = (x+1)r_{n-1}(x,-1) + \left((-1)^{n-1}-1\right)xr_{n-2}(x,-1)$ with initial values $r_0(x,-1)=1$ and $r_1(x,-1) = 1+x$. This gives by induction

$$\begin{aligned} r_{2n}(x,-1) &= \left(1+x^2\right)^n, \\ r_{2n+1}(x,-1) &= (1+x)\left(1+x^2\right)^n. \end{aligned} \tag{4.12}$$

Of course this also follows from (3.6).

In the general case we get from the defining relation $\begin{bmatrix} n \\ k \end{bmatrix} = q^k \begin{bmatrix} n-1 \\ k \end{bmatrix} + \begin{bmatrix} n-1 \\ k-1 \end{bmatrix}$ of the $q-$ Pascal matrix

**Proposition 4.3**

The equivalence classes $\varepsilon(n,k) = e(n,k) + qo(n,k)$ satisfy
$$\varepsilon(n,k) \equiv q^k\left(\varepsilon(n-1,k) + \varepsilon(n-1,k-1)\right) \tag{4.13}$$
with $\varepsilon(n,0) = 1$ for all $n$ and $\varepsilon(0,k) = 0$ for $k > 0$
and the equivalence classes $\lambda(n,k) = L(n,k) + q\overline{L}(n,k)$ satisfy
$$\lambda(n,k) \equiv q^k\lambda(n-1,k) + \lambda(n-1,k-1) \tag{4.14}$$
with $\lambda(n,0) = 1$ for all $n$ and $\lambda(0,k) = 0$ for $k > 0$.

**Remark**

It is clear that $\lambda(n,k)$ and $\varepsilon(n,k)$ are uniquely determined by these rules. Thus we get a very simple approach to the matrices $(L(n,k))$ and $(e(n,k))$ and their properties.

From
$\varepsilon(n,k) \equiv q^k\varepsilon(n-1,k) + q^k\varepsilon(n-1,k-1) \equiv \varepsilon(n-2,k) + \varepsilon(n-2,k-1) + q\left(\varepsilon(n-2,k-1) + \varepsilon(n-2,k-2)\right)$
$\equiv \varepsilon(n-2,k) + \binom{n-2}{k-1} + q\varepsilon(n-2,k-2)$

we get $e(n,k) = e(n-2,k) + \binom{n-2}{k-1} + \binom{n-2}{k-2} - e(n-2,k-2)$

and thus again (2.1).
In the same way we get (3.4).

Formula (4.14) implies that

$$L(n, 2k) = L(n-1, 2k) + L(n-1, 2k-1),$$
$$L(n, 2k+1) = \overline{L}(n-1, 2k+1) + L(n-1, 2k). \qquad (4.15)$$

This is obvious for the interpretation of $L(n,k)$ as the set of $w \in W_{n,k}$ with an even number of inversions if we consider the last entry of $w$.

**4.4**
Let us give some more examples:

**a)** The recurrence (4.7) implies
$$L_{2n+1}(x) + q\overline{L}_{2n+1}(x) \equiv (x+1)\left(L_{2n}(x) + q\overline{L}_{2n}(x)\right)$$
which again gives (3.7)
and
$$L_{2n}(x) + q\overline{L}_{2n}(x) \equiv (x+1)\left(L_{2n-1}(x) + q\overline{L}_{2n-1}(x)\right) + (q-1)x\left(L_{2n-2}(x) + q\overline{L}_{2n-2}(x)\right)$$
which implies

$$L_{2n}(x) = (x+1)L_{2n-1}(x) + x\overline{L}_{2n-2}(x) - xL_{2n-2}(x) = (x+1)L_{2n-1}(x) - 2xL_{2n-2}(x) + x(1+x)^{2n-2}.$$

**b)** The well-known identity $\sum_{k=0}^{n} q^{k^2} \begin{bmatrix} n \\ k \end{bmatrix}^2 = \begin{bmatrix} 2n \\ n \end{bmatrix}$ implies

$$\sum_{k=0}^{n} q^{k^2} \begin{bmatrix} n \\ k \end{bmatrix}^2 \equiv \sum_{k} \begin{bmatrix} n \\ 2k \end{bmatrix}^2 + q \sum_{k} \begin{bmatrix} n \\ 2k+1 \end{bmatrix}^2$$
$$\equiv \sum_{k} \left(L(n, 2k) + q\overline{L}(n, 2k)\right)^2 + q \sum_{k} \left(L(n, 2k+1) + q\overline{L}(n, 2k+1)\right)^2 \equiv L(2n, n) + q\overline{L}(2n, n)$$
which gives
$$\sum_{k} \left(L(n, 2k)^2 + \overline{L}(n, 2k)^2 + 2L(n, 2k+1)\overline{L}(n, 2k+1)\right) = L(2n, n). \qquad (4.16)$$

**c)** A more transparent proof of (3.12) runs as follows: The well-known identity
$$\sum_{n \geq 0} \begin{bmatrix} n \\ k \end{bmatrix} x^n = \frac{x^k}{(1-x)(1-qx)\cdots(1-q^k x)} \text{ implies}$$

$$(1+qx)\cdots(1+q^k x)(1-x)(1-qx)\cdots(1-q^k x)\sum_{n \geq 0}\begin{bmatrix}n\\k\end{bmatrix}x^n = (1-x)(1-q^2 x^2)\cdots(1-q^{2k}x^2)\sum_{n \geq 0}\begin{bmatrix}n\\k\end{bmatrix}x^n$$

$$= x^k \sum_{k=0}^{n} q^{\binom{k+1}{2}} \begin{bmatrix} n \\ k \end{bmatrix} x^k.$$

If we reduce this identity modulo $q^2 - 1$ we get (3.12).

The same proof also gives again (2.13). We have only to observe the parity of $\binom{4k+j}{2}$.

**d)** The $q$-Fibonacci polynomials $F_n(s,q) = \sum_{k=0}^{\lfloor \frac{n-1}{2} \rfloor} q^{k(k-1)} \begin{bmatrix} n-1-k \\ k \end{bmatrix} s^k$

satisfy $F_n(s,q) = F_{n-1}(s,q) + q^{n-3} s F_{n-2}(s,q)$ (cf. [3]).

For $q = 1$ we get the Fibonacci polynomials $F_n(s) = \sum_{k=0}^{\lfloor \frac{n-1}{2} \rfloor} \binom{n-1-k}{k} s^k$.

Note that $(F_n(1))_{n \geq 0} = (0, 1, 1, 2, 3, 5, 8, 13, 21, \cdots)$ gives the Fibonacci numbers.

For $q = -1$ we get

$$F_{2n}(s,-1) = F_n(s^2),$$
$$F_{2n+1}(s,-1) = F_{n+1}(s^2) + s F_n(s^2).$$
(4.17)

Therefore the polynomials

$$\varphi_n(s,q) = \sum_{k=0}^{\lfloor \frac{n-1}{2} \rfloor} \lambda(n-1-k, k) s^k = F_n(s,q) \bmod (q^2 - 1) \text{ satisfy}$$

$$\varphi_n(s,q) \equiv \varphi_{n-1}(s,q) + q^{n-3} s \varphi_{n-2}(s,q) \bmod (q^2 - 1).$$

The first terms are

$0, 1, 1, 1+s, 1+s+sq, 1+2s+s^2+sq, 1+2s+2s^2+(2s+s^2)q, \cdots.$

Let now $\varphi_n(s,q) = f_n(s) + q \overline{f}_n(s)$. Then

$$f_n(s) = \sum_{k=0}^{\lfloor \frac{n-1}{2} \rfloor} L(n-1-k, k) s^k.$$
(4.18)

We get

$f_{2n+1}(s) = f_{2n}(s) + s f_{2n-1}(s),$
$f_{2n}(s) = f_{2n-1}(s) + s \overline{f}_{2n-2}(s).$

This gives by induction or by considering (4.17)

$$f_{2n}(s) = \frac{F_{2n}(s) + F_n(s^2)}{2},$$
$$f_{2n+1}(s) = \frac{F_{2n+1}(s) + F_{n+1}(s^2) + s F_n(s^2)}{2}.$$
(4.19)

For $s=1$ we get an analogue of the Fibonacci numbers

$$f_n(1) = \sum_{k=0}^{\left\lfloor \frac{n-1}{2} \right\rfloor} L(n-1-k,k). \qquad (4.20)$$

The first terms are

$$(f_n(1))_{n\geq 0} = (0,1,1,2,2,4,5,9,12,\cdots).$$

This sequence also occurs in OEIS[7], A102526.

From (4.19) we get $f_{2n}(1) = \dfrac{F_{2n}+F_n}{2}$, $f_{2n+1}(1) = \dfrac{F_{2n+1}+F_{n+2}}{2}$.

**e)** As is well known the sequence

$$(F_n(-1))_{n\geq 0} = \left( \sum_{k=0}^{\left\lfloor \frac{n-1}{2} \right\rfloor} (-1)^k \binom{n-1-k}{k} \right)_{n\geq 0} = (0,1,1,0,-1,-1,0,1,1,0,-1,-1,\cdots) \text{ is periodic with}$$

period 6.

Consider the $q-$analogue $f(n) = \sum_{k=0}^{n} (-1)^k q^{\binom{k+1}{2}} \begin{bmatrix} n-1-k \\ k \end{bmatrix}$. The first terms are

$\{0, 1, 1, 1-q, 1-q-q^2, 1-q-q^2, 1-q-q^2+q^5, 1-q-q^2+q^5+q^7, 1-q-q^2+q^5+q^7, 1-q-q^2+q^5+q^7-q^{12},$
$1-q-q^2+q^5+q^7-q^{12}-q^{15}, 1-q-q^2+q^5+q^7-q^{12}-q^{15}, 1-q-q^2+q^5+q^7-q^{12}-q^{15}+q^{22}\}$

Note that these are partial sums of Euler's pentagonal number series

$$\prod_{n\geq 1}(1-q^n) = 1-q-q^2+q^5+q^7-q^{12}-q^{15}+q^{22}+q^{26}--++\cdots.$$

As shown in [4] it satisfies $f(n) = f(n-1) - q^{n-2}f(n-3) + q^{n-2}f(n-4)$.
Let $\varphi(n) = f(n) \bmod(q^2-1)$.
Then it is easy to verify that $\varphi(n+12) = \varphi(n)$ with initial values
$0,1,1,1-q,-q,-q,0,q,q,q,-1,-1,-1$.

Therefore we get

**Proposition 4.4**

*The sequence* $\left( \sum_{k=0}^{\left\lfloor \frac{n-1}{2} \right\rfloor} (-1)^k e(n-1-k,k) \right)_{n\geq 0}$ *is periodic with period* 12 *with initial values*

$0,1,1,1,0,0,0,0,0,-1,-1,-1$.

## 5. A more general case

### 5.1

Let $e(n,k,j,p)$ be the number of $k-$ subsets of $\{1,2,\cdots,n\}$ whose sums are congruent to $j$ modulo $p$.

Then

$$\prod_{j=1}^{n}(1+q^j x) = \sum_k x^k \sum_{j_1<j_2<\cdots<j_k} q^{j_1+\cdots+j_k} \equiv \sum_k x^k \sum_{j=0}^{p-1} e(n,k,j,p) q^j \mod(q^p-1).$$

On the other hand we know that

$$\prod_{j=1}^{n}(1+q^j x) = \sum_{k=0}^{n} q^{\binom{k+1}{2}} \begin{bmatrix} n \\ k \end{bmatrix}_q x^k.$$

Therefore we get

$$e(n,k,j,p) = [q^j] q^{\binom{k+1}{2}} \begin{bmatrix} n \\ k \end{bmatrix}_q \mod(q^p-1). \qquad (4.21)$$

Let now

$$\varepsilon(n,k,p,q) \equiv \sum_{j=0}^{p-1} e(n,k,j,p) q^j \equiv q^{\binom{k+1}{2}} \begin{bmatrix} n \\ k \end{bmatrix}_q \mod(q^p-1). \qquad (4.22)$$

The equivalence classes $\varepsilon(n,k,p,q)$ satisfy

$$\varepsilon(n,k,p,q) \equiv q^k \left( \varepsilon(n-1,k,p,q) + \varepsilon(n-1,k-1,p,q) \right) \mod(q^p-1) \qquad (4.23)$$

with $\varepsilon(n,0,p,q)=1$ for all $n$ and $\varepsilon(0,k,p,q)=0$ for $k>0$.

For example the first terms of $\left(\varepsilon(n,k,3,q)\right)$ are

$$\begin{pmatrix}
1 & 0 & 0 & 0 & 0 & 0 & 0 \\
1 & q & 0 & 0 & 0 & 0 & 0 \\
1 & q+q^2 & 1 & 0 & 0 & 0 & 0 \\
1 & 1+q+q^2 & 1+q+q^2 & 1 & 0 & 0 & 0 \\
1 & 1+2q+q^2 & 2+2q+2q^2 & 2+q+q^2 & q & 0 & 0 \\
1 & 1+2q+2q^2 & 4+3q+3q^2 & 4+3q+3q^2 & 1+2q+2q^2 & 1 & 0 \\
1 & 2+2q+2q^2 & 5+5q+5q^2 & 8+6q+6q^2 & 5+5q+5q^2 & 2+2q+2q^2 & 1
\end{pmatrix}$$

## 5.2

Let us first derive an analogue of (2.2):

**Proposition 5.1**

Let $e_n(j,x,p) = \sum_{k=0}^{n} e(n,k,j,p)x^k$. For $n \geq p$ and each $j$ we get

$$e_n(j,x,p) = (1+x^p)e_{n-p}(j,x,p) + \frac{1}{p}\sum_{i=1}^{p-1}\binom{p}{i}x^i(1+x)^{n-p}. \tag{4.24}$$

**Proof**

To show this observe first that $e(p,k,j,p) = e(p,k,0,p)$ for all $j$ with $0 \leq j \leq p$ if $1 \leq k \leq p-1$. This implies $e(p,k,j,p) = \frac{1}{p}\binom{p}{k}$ for all $j$.

J. Schoissengeier has shown me a simple proof due to V. Losert: Consider the group $\mathbb{Z}_p$ of integers modulo $p$ and let $E(p,k,i)$ be the set of all $k$-subsets $S$ of $\mathbb{Z}_p$ whose sum of elements is congruent to $i$ modulo $p$. Let $S+1$ be the the set of elements $s+1$ with $s \in S$. The map $S \to S+1$ gives a bijection from $E(p,k,i)$ to $E(p,k,i+k)$. Iterating this map we see that all sets $E(p,k,i)$ have the same size because the multiples of $k$ exhaust $\mathbb{Z}_p$.

Now consider the elements $i$ with $n-p+1 \leq i \leq n$.
The number of $k$-sets whose sums are congruent to $j$ modulo $p$ which contain none of these numbers is $e(n-p,k,j,p)$.
Since the sum of $p$ consecutive numbers is divisible by $p$ the number of $k$-sets which contain all elements $n-p+1 \leq i \leq n$ is $e(n-p,k-p,j,p)$.
Let now $1 \leq \ell \leq p-1$.
For each $(k-\ell)$-subset of $\{1,2,\cdots,n-p\}$ there are precisely $\frac{1}{p}\binom{p}{\ell}$ $\ell$-subsets of $\{n-p+1,\cdots n\}$ such hat their union is $\equiv j \bmod p$ because $p$ consecutive numbers contain all residues modulo $p$. This gives $\binom{n-p}{k-\ell}\frac{1}{p}\binom{p}{\ell}$ sets which are $\equiv j \bmod p$.

This gives for $n \geq p$

$$e(n,k,j,p) = e(n-p,k,j,p) + \frac{1}{p}\sum_{\ell=1}^{p-1}\binom{p}{\ell}\binom{n-p}{k-\ell} + e(n-p,k-p,j,p),$$

which implies (4.24).

By applying (4.24) to $n$ and $n-1$ we get the homogeneous recursion

$$e_n(j,x,p) - (1+x)e_{n-1}(j,x,p) - (1+x^p)e_{n-p}(j,x,p) + (1+x)(1+x^p)e_{n-p-1}(j,x,p) = 0. \quad (4.25)$$

By summing over all $j$ we get

$$\varepsilon(n,k,p,q) = \varepsilon(n-p,k,p,q) + \frac{1}{p}\sum_{\ell=1}^{p-1}\binom{p}{\ell}\binom{n-p}{k-\ell}(1+q+\cdots+q^{p-1}) + \varepsilon(n-p,k-p,p,q). \quad (4.26)$$

For a primitive $p$-th root of unity $q = \zeta$ (4.26) reduces to

$$\varepsilon(n,k,p,\zeta) = \varepsilon(n-p,k,p,\zeta) + \varepsilon(n-p,k-p,p,\zeta). \quad (4.27)$$

### 5.3

Let us now derive an analogue of (2.4).

Let $\zeta$ be a primitive $p$-th root of unity. Then

$$\prod_{j=1}^{n}(1+\zeta^j x) = \sum_k x^k \sum_{j=0}^{p-1} e(n,k,j,p)\zeta^j.$$

Observe that

$$\sum_{\ell=0}^{p-1}\prod_{j=1}^{n}(1+\zeta^{\ell j}x) = \sum_k x^k \sum_{\ell=0}^{p-1}\sum_{j=0}^{p-1} e(n,k,j,p)\zeta^{\ell j} = p\sum_k e(n,k,0,p)x^k.$$

On the other hand we have

$$\sum_{\ell=0}^{p-1}\prod_{j=1}^{n}(1+\zeta^{\ell j}x) = (1+x)^n + \sum_{\ell=1}^{p-1}\prod_{j=1}^{n}(1+\zeta^{\ell j}x).$$

Since each product of $1+\zeta^{\ell j}x$ over $p$ consecutive values of $j$ equals $1+x^p$ we see that

$$\sum_{\ell=1}^{p-1}\prod_{j=1}^{pn+i}(1+\zeta^{\ell j}x) = b_i(x)(1+x^p)^n$$

for some polynomial $b_i(x)$ of degree $i$.

Therefore the polynomial $e_n(x,p) = \sum_k e(n,k,0,p)x^k$ satisfies

$$e_{pn+i}(x,p) = \frac{(1+x)^{pn+i} + b_i(x)(1+x^p)^n}{p}. \quad (4.28)$$

This implies

$$\sum_{n\geq 0} e_n(x,p)z^n = \frac{1}{p}\left( \frac{1}{1-(x+1)z} + \frac{\sum_{i=0}^{p-1} b_i(x)z^i}{1-(1+x^p)z^p} \right). \quad (4.29)$$

Observe that this also follows from (4.25).

## 5.4

We define now $p$ – **Losanitsch numbers** $L(n,k,j,p)$ by

$$L(n,k,j,p) = [q^j]\begin{bmatrix} n \\ k \end{bmatrix}_q \mod (q^p - 1) \quad (4.30)$$

or equivalently by

$$\lambda(n,k,p,q) \equiv \sum_{j=0}^{p-1} L(n,k,j,p)q^j \equiv \begin{bmatrix} n \\ k \end{bmatrix}_q \mod (q^p - 1). \quad (4.31)$$

The equivalence classes $\lambda(n,k,p,q)$ satisfy

$$\lambda(n,k,p,q) \equiv q^k \lambda(n-1,k,p,q) + \lambda(n-1,k-1,p,q) \mod (q^p - 1). \quad (4.32)$$

with $\lambda(n,0,p,q) = 1$ for all $n$ and $\lambda(0,k,p,q) = 0$ for $k > 0$.

For example the first terms of $(\lambda(n,k,3))$ are

$$\begin{pmatrix}
1 & 0 & 0 & 0 & 0 & 0 & 0 \\
1 & 1 & 0 & 0 & 0 & 0 & 0 \\
1 & 1+q & 1 & 0 & 0 & 0 & 0 \\
1 & 1+q+q^2 & 1+q+q^2 & 1 & 0 & 0 & 0 \\
1 & 2+q+q^2 & 2+2q+2q^2 & 2+q+q^2 & 1 & 0 & 0 \\
1 & 2+2q+q^2 & 4+3q+3q^2 & 4+3q+3q^2 & 2+2q+q^2 & 1 & 0 \\
1 & 2+2q+2q^2 & 5+5q+5q^2 & 8+6q+6q^2 & 5+5q+5q^2 & 2+2q+2q^2 & 1
\end{pmatrix}$$

The polynomials $L_n(x,p,j) = \sum_{k=0}^{n} L(n,k,j,p)x^k$ will be called $p$ – Losanitsch polynomials.
Note that all polynomials $L_n(x,p,j)$ are palindromic.

Since $e(n,k,j,p) = [q^j]q^{\binom{k+1}{2}}\begin{bmatrix} n \\ k \end{bmatrix} \mod (q^p - 1)$ we see that

$$L(n,k,j,p) \equiv e\left(n,k,j+\binom{k+1}{2},p\right) \mod (q^p - 1).$$

Therefore each column $(L(n,k,j,p))_{n\geq 0}$ coincides with a column $(a(n,k,i,p))_{n\geq 0}$ for some $i$.

Therefore we get

**Proposition 5.2**

*The matrix* $\left(L(n,k,0,p)\right)_{n,k\geq 0}$ *is the uniquely determined matrix whose columns are* $\left(e(n,k,i,p)\right)_{n\geq 0}$ *for some* $i$ *and where each element of the main diagonal is* 1.

From (4.31) and (4.8) we conclude

**Theorem 5.3**

$L(n,k,0,p)$ *is the number of elements* $w \in W_{n,k}$ *such that* $\mathrm{inv}(w)$ *is a multiple of* $p$.